\magnification=\magstep1
\input amssym.def
\input amssym.tex
\input psfig.sty
\overfullrule=0pt 
\font\titlefont=cmssbx10 scaled \magstep2
\font\caps=cmcsc10	
\def\balpha{\mathop{\alpha\kern-.6em{\alpha}}\nolimits}		

\def\ref#1{{\rm [#1]}} 
\def\Refs{\bigskip\baselineskip=12pt\frenchspacing 
	\leftline{\bf References}\smallskip} 
 
\def\eop{\hbox{\vrule width6pt height7pt depth1pt}}
\def\qed{~\hfill~\eop\medskip} 
\def\complex{{\Bbb C}}
\def\FF{{\Bbb F}} 
\def\que{{\Bbb Q}}
\def\real{{\Bbb R}} 
\def\P{{\cal P}}
\def\L{{\cal L}}
\def\ov{\overline}
\def\varep{\varepsilon}
\def\iitem{\itemitem}
\def\dist{\mathop{\rm dist}\nolimits}
\def\frac#1#2{{\textstyle{#1\over#2}}}

\outer\def\demo #1. #2\par{\medbreak\noindent {\bf#1.\enspace}
	{\rm#2}\par\ifdim\lastskip<\medskipamount
	\removelastskip\penalty55\medskip\fi}
\def\newpage{\vfill\eject}
\topinsert\vskip.5in\endinsert
\centerline{\titlefont Totally real origami and}
\centerline{\titlefont impossible paper folding}
\bigskip
\centerline{\caps David Auckly and John Cleveland}
\smallskip
\centerline{Department of Mathematics}
\centerline{The University of Texas at Austin}
\centerline{Austin, TX 78712}
\vskip.3in 

\baselineskip=17pt 

Origami is the ancient Japanese art of paper folding. 
It is possible to fold many intriguing geometrical shapes with paper \ref{M}. 
In this article, the question we will answer is which shapes are possible 
to construct and which shapes are impossible to construct using origami. 
One of the most interesting things we discovered is that it is impossible 
to construct a cube with twice the volume of a given cube using origami, 
just as it is impossible to do using a compass and straight edge. 
As an unexpected surprise, our algebraic characterization of origami is 
related to David Hilbert's $17^{th}$ problem. Hilbert's 
problem is to show that any 
rational function which is always non-negative is a sum of squares of 
rational functions \ref{B}. This problem was solved by Artin in 1926 
\ref{Ar}. We would like to thank John Tate for noticing the relationship 
between our present work and Hilbert's $17^{th}$ problem. 
This research is the result of a project in the Junior Fellows Program 
at The University of Texas. The Junior Fellows Program is a program in 
which a junior undergraduate strives to do original research under the 
guidance of a faculty mentor. 

The referee mentioned two references which the reader may find 
interesting. ``Geometric Exercises in Paper Folding'' addresses practical 
problems of paper folding \ref{R}. Among many other things, Sundara Row 
gives constructions for the 5-gon, the 17-gon, and duplicating a cube. 
His constructions, however, use more general folding techniques than 
the ones we consider here. Felix Klein cites Row's work in his lectures on 
selected questions in elementary geometry \ref{K}. 

In order to understand the rules of origami construction, we will first 
consider a sheet of everyday notebook paper. Our work with notebook paper 
will serve as an intuitive model for our definition of origami constructions 
in the Euclidean plane. There are four natural methods of folding a piece 
of paper. The methods will serve as the basis of the definition of an 
origami pair. 


\hskip 90bp\psfig{file=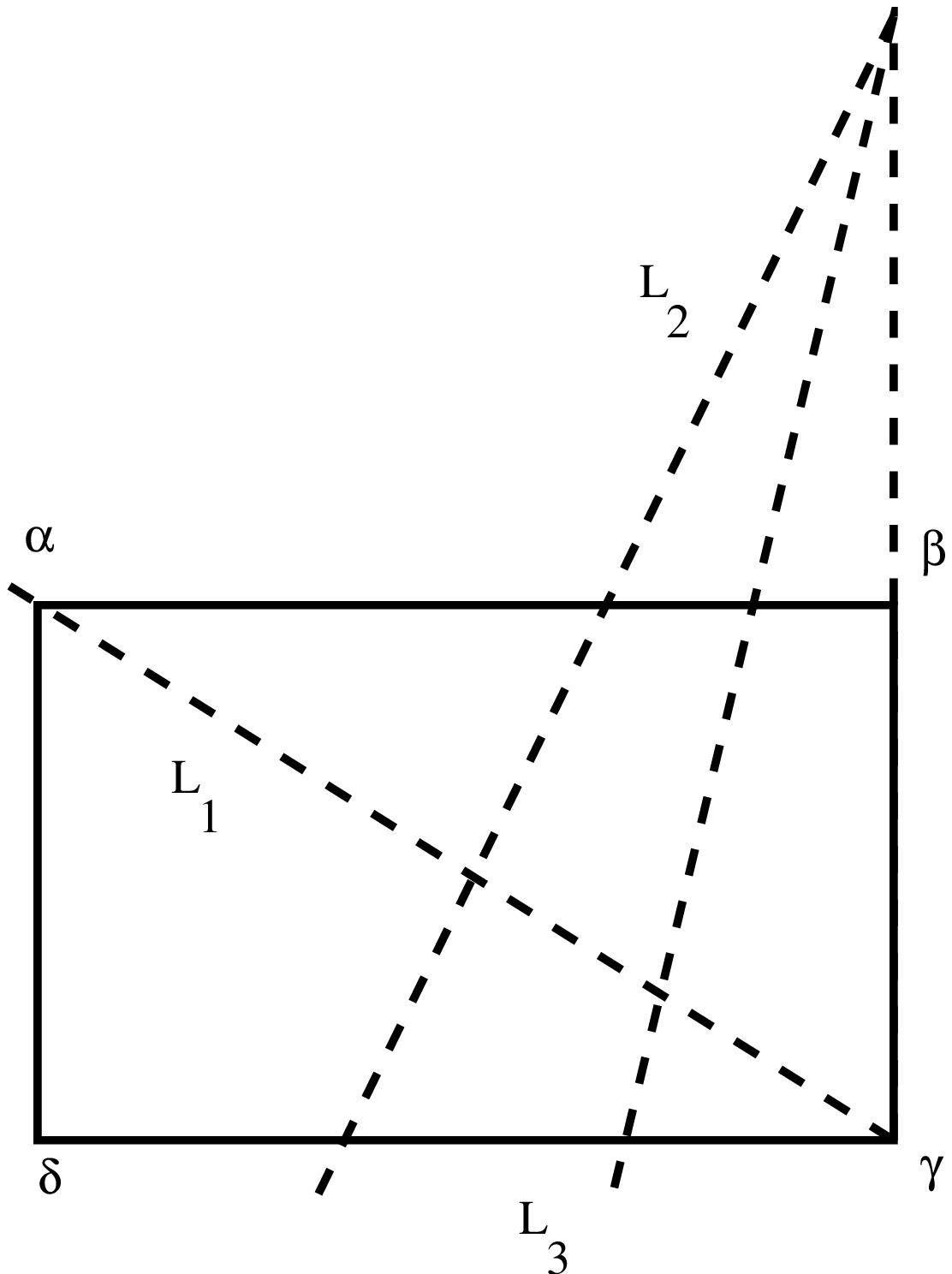,width=2.8truein}

\medskip
\centerline{Figure 1}
\smallskip


We construct the line $L_1$, by folding a crease between two 
different corners of the paper. Another line may be constructed by 
matching two corners. For example, if corners $\alpha$ and $\gamma$ are 
matched, the crease formed, $L_2$, will be the perpendicular bisector 
of the segment $\overline{\alpha\gamma}$. another natural 
construction is matching one line to another line. 
For instance, $\overline{\beta\gamma}$, the paper's edge, and $L_2$ 
are lines. 
If we lay $\overline{\beta\gamma}$ upon $L_2$ and form the crease, then 
we obtain $L_3$ which is the angle bisector of the two lines. 
If we start with two parallel lines in this third construction, then 
we will just get a parallel line half way in between. 

The fourth and final construction which seems natural is consecutive 
folding. This is similar to rolling up the sheet of paper only one 
does not roll it up,  he folds it up. 
More explicitly, start with a piece of paper with two creases on it as in 
figure~2. 
Fold along line $L_1$ and do not unfold the piece of paper. 
Notice that line $L_2$ lies over the sheet of paper. 
With the paper still folded, fold the sheet of paper along the crease 
$L_2$ to obtain a new crease on the sheet underneath $L_2$. 
If we name this new crease $L_3$ and unfold the sheet of paper, 
then it is easy to see that line $L_3$ is the mirror image or reflection 
of line $L_2$ about line $L_1$. 
\newpage 


\hskip 80bp\psfig{file=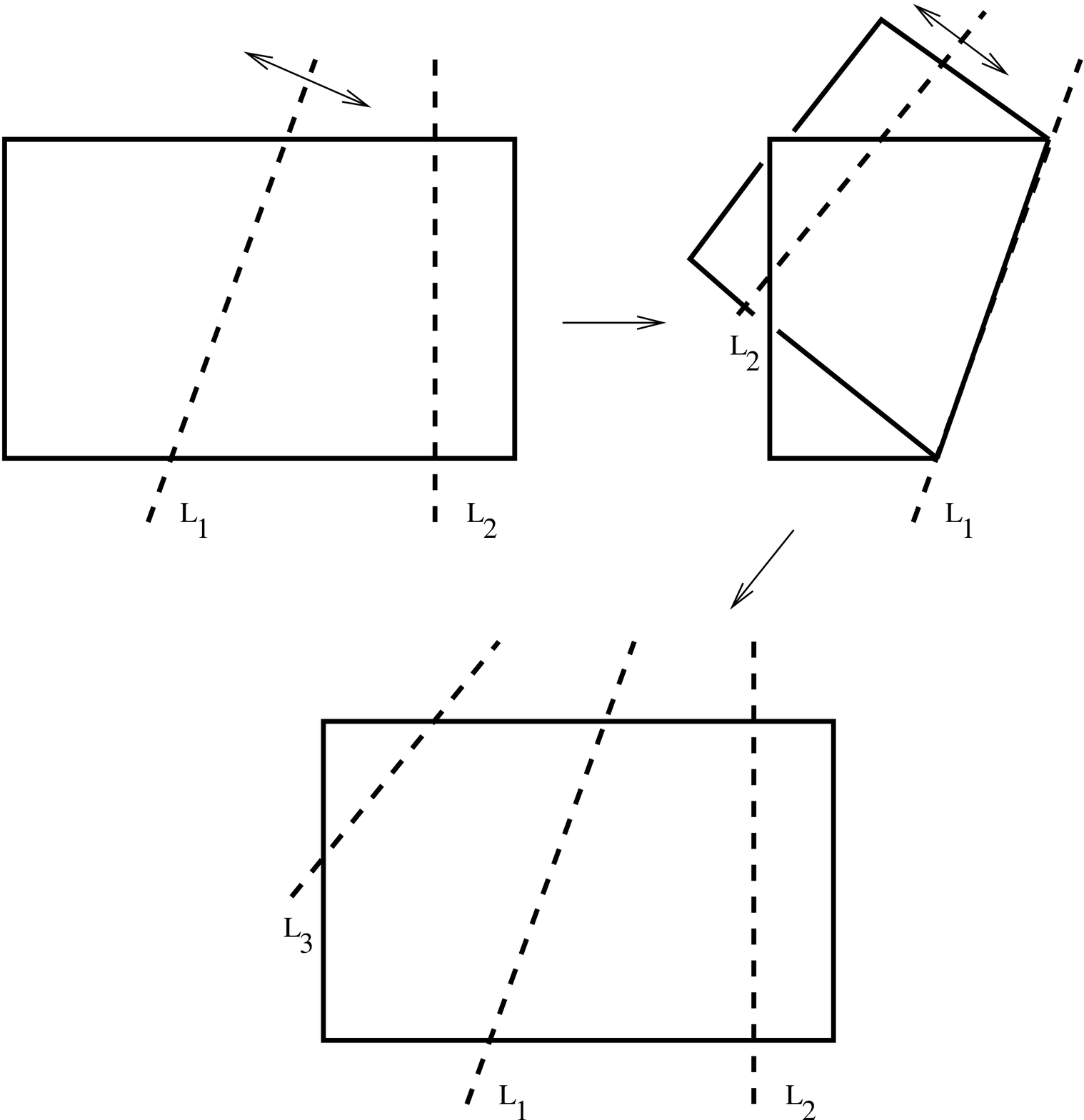,width=3.4truein}

\medskip
\centerline{Figure 2}
\smallskip

We now formalize these methods to define an origami pair on the plane. 
The creases on our sheet of paper are merely lines in the plane, and 
the corners of the paper are represented by points where lines (creases) meet. 
This previous discussion is the motivation for the following 
definition.

\demo Definition. 
$(\P,\L)$ is an {\it origami pair\/} if $\P$ is a set of points in $\real^2$ 
and $\L$ is a collection of lines in $\real^2$ satisfying: 
\smallskip
\iitem{i)}  The point of intersection of any two non-parallel lines 
in $\L$ is a point in $\P$. 
\iitem{ii)} Given any two  distinct points in $\P$, there is a line in $\L$ 
going through them. 
\iitem{iii)} Given any two distinct points in $\P$, the perpendicular 
bisector of the line segment with given end points is a line in $\L$. 
\iitem{iv)} If $L_1$ and $L_2$ are lines in $\L$, then the line which is 
equidistant from $L_1$ and $L_2$ is in $\L$. 
\iitem{v)} If $L_1$ and $L_2$ are lines in $\L$, then there exists a line 
$L_3$ in $\L$ such that $L_3$ is the mirror reflection of $L_2$ about $L_1$. 
\smallskip

For any subset of the plane containing at least two points, there is at 
most one collection of lines which will pair with it to become an 
origami pair. 

\demo Definition. 
A subset of $\real^2,\P$, is {\it closed under origami constructions\/} if 
there exists a collection of lines, $\L$, such that $(\P,\L)$ is an 
origami pair. 

The question which we answer in this paper is which points may be constructed 
from just two points, using only the origami constructions described above. 
We will call that collection of points the set of origami constructible points. 

\demo Definition. 
$\P_0 = \cap \{\P\mid (0,0), (0,1) \in \P$ and  $\P$ is closed under 
origami constructions$\}$ is the set of {\it origami constructible points\/}. 

Before we explain the structure of $\P_0$, we give an example of an origami 
construction analogous to many compass and straight edge constructions, 
namely, the construction of parallel lines. 

\proclaim Lemma. 
It is possible to construct a line parallel to a given line through 
any given point using origami. 


\hskip 110bp\psfig{file=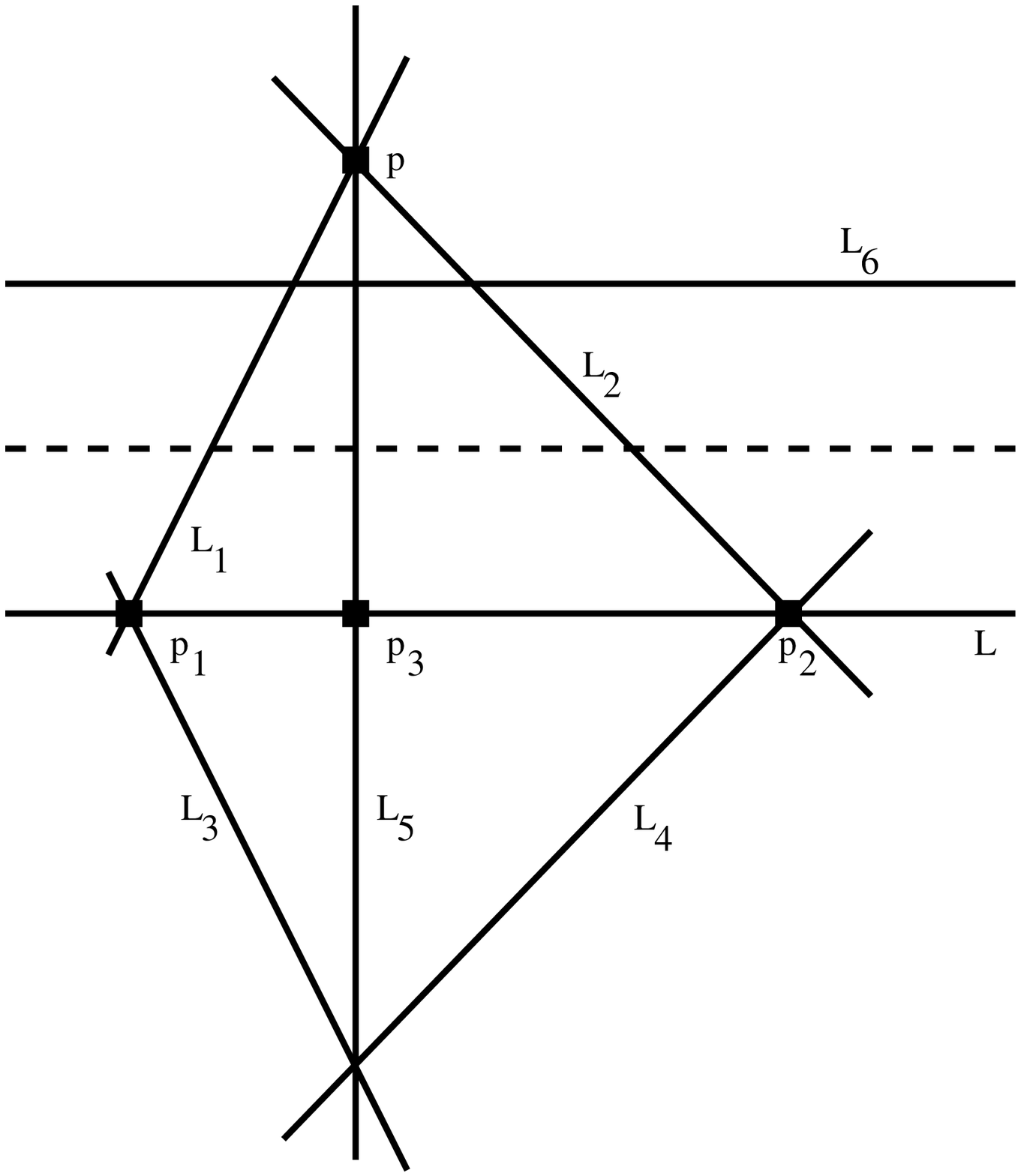,width=2.5truein}

\medskip
\centerline{Figure 3}
\smallskip

\demo Proof. 
Refer to Figure 3. Given a line $L$ and a point $p$, pick two points 
$p_1$ and $p_2$ on $L$. 
By property~ii) in the definition of an origami pair, we may construct 
lines $L_1$ and $L_2$ running through $p_1,p$ and $p_2,p$,  respectively. 
By property~v) we may reflect $L_1$ and $L_2$ through $L$ to obtain 
$L_3$ and $L_4$. 
Now the intersection of $L_3$ and $L_4$ is a constructible point, so 
there is a line, $L_5$ through this point and the given point, $p$, 
by properties~i) and ii). 
Call the point where $L_5$ and $L$ intersect $p_3$. 
To finish the construction, use property~iii) to construct a perpendicular 
bisector to $\overline{p,p_3}$, and reflect $L$ through this bisector with 
property~v) to obtain the desired line, $L_6$. 
It is a straightforward exercise to show that $L_6$ has the desired 
properties.\qed  

The reader may wish to try some constructions on his own. Two especially 
interesting exercises to attempt are the construction of a right triangle 
with given legs and the construction  of a right triangle with a given 
hypotenuse and leg. 
More explicitly, 
given four distinct points $\alpha,\beta,\gamma$ and $\delta$, the reader may 
try to construct a point $\varep$ such that $\alpha,\beta,\varep$ are 
the vertices of a right triangle with legs $\ov{\alpha\beta}$ and 
$\ov{\beta\varep}$ such that the length of $\ov{\beta\varep}$  equals the 
length of $\ov{\gamma\delta}$. 

Now that we have a better feel for origami constructions, we will start 
developing tools  to show that some figures are not constructible. 
The first thing we need is the notion of an origami number. 

\demo Definition. 
$\FF_0 = \{\alpha\in\real\mid \exists\ v_1,v_2\in \P$  such that 
$|\alpha| = \dist (v_1,v_2)\}$ is the {\it set of origami numbers\/}. 

It is easy to see that $(x,y)\in \P_0$ if and only if $x$ and $y$ are 
both in $\FF_0$. It is also easy to see that the numbers $\frac12,\frac14,
\frac18,\ldots$ are origami numbers. 
To see that $\frac15$ is an origami number consider figure~4. 

\hskip 100bp\psfig{file=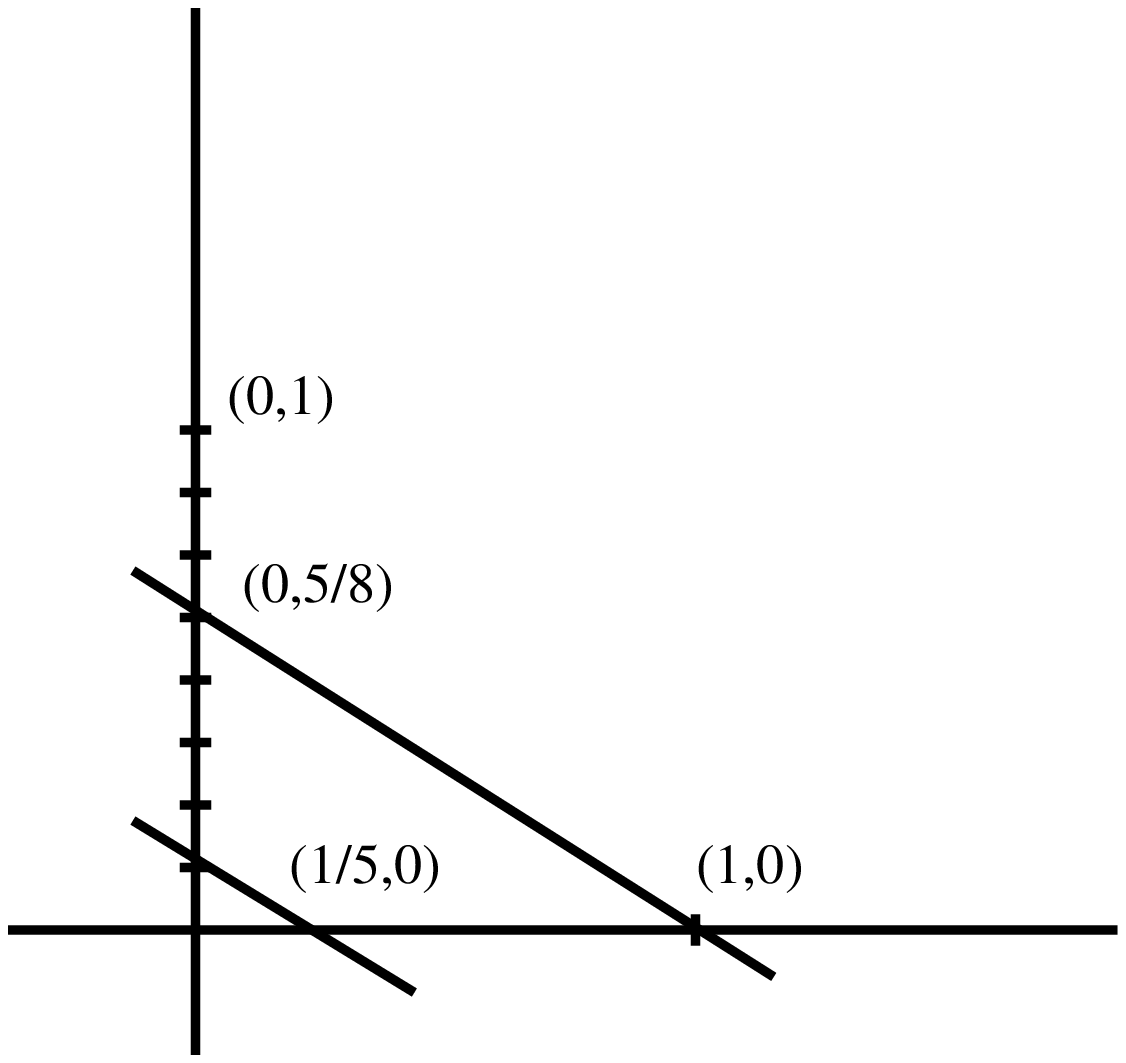,width=2.2truein}

\medskip
\centerline{Figure 4}
\smallskip 
In figure 4 a line through $(0,\frac58)$ and $(1,0)$ is constructed, then 
a parallel line through $(0,\frac18)$ is constructed. 
This parallel line intersects the $x$-axis at $(\frac15,0)$, therefore 
$\frac15$ is an origami number. 
Another class of origami numbers can be generalized by a simple 
geometric construction. 
Starting with any segment, it is possible to construct a right triangle 
as in figure~5. 

\hskip 90bp\psfig{file=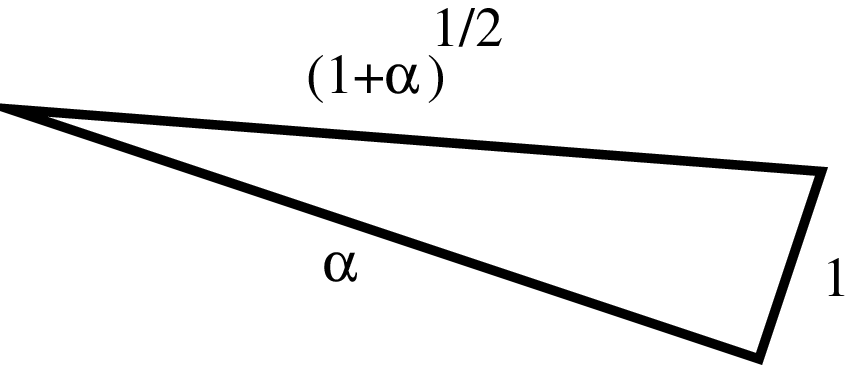,width=2.8truein}

\medskip
\centerline{Figure 5}
\smallskip
\noindent
It follows that $\sqrt{1+\alpha^2}$ is an origami number whenever $\alpha$ 
is an origami number. Using this construction, we see that 
$$\sqrt2 = \sqrt{1+1^2}\quad\hbox{and}\quad 
\sqrt3 = \sqrt{1+(\sqrt2\,)^2}$$ 
are origami numbers. 
In fact, the sum, difference, product, and quotient of origami numbers 
are origami numbers. 

\proclaim Theorem. 
The collection of origami numbers, $\FF_0$ is a field closed under 
the operation $\alpha\mapsto \sqrt{1+\alpha^2}$. 

\demo Proof. 
If $\alpha,\beta\in \FF_0$, it follows from the definition that 
$-\alpha\in\FF_0$ and it is easy to show that $\alpha+\beta\in\FF_0$. 
Straightforward constructions with similar triangles are enough to show
that $\alpha\cdot\beta$, $\alpha^{-1}\in\FF_0$. See figure~6. 

\hskip 40bp\psfig{file=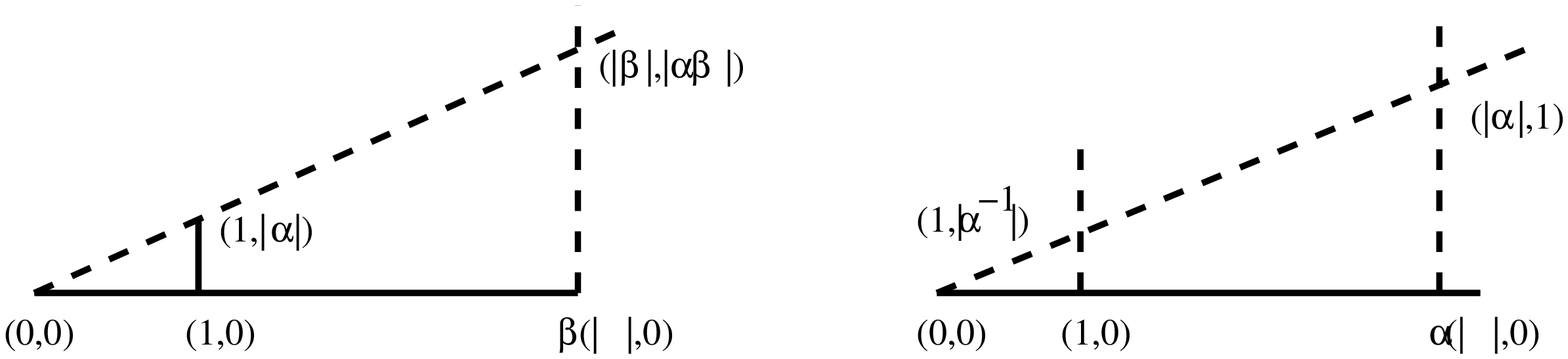,width=5.2truein}

\medskip
\centerline{Figure 6}
\smallskip
\noindent
In the discussion preceding this theorem we showed that $\sqrt{1+\alpha^2}$ 
is an origami number whenever $\alpha$ is. 
The proof is therefore complete.\qed 

Now that we have some algebraic operations which will produce origami numbers, 
it is natural to ask if there are any more operations which will produce 
origami numbers. 
Once we have a list of all ways to create origami numbers and a method 
to test if a given number can be achieved, then we will know which geometric 
shapes are constructible, and which shapes are not constructible. 
This is because any figure is constructible if and only if the coordinates 
of all of the vertices are origami numbers. 

\demo Definition. 
$\FF_{\sqrt{1+x^2}}$ is the smallest subfield of $\complex$ closed under the 
operation $x\mapsto \sqrt{1+x^2}$. 

The preceding Theorem may be rephrased as $\FF_{\sqrt{1+x^2}} \subset \FF_0$. 
It is in fact true that $\FF_0 = \FF_{\sqrt{1+x^2}}$. 
Thus, the previously listed operations which produce origami numbers are 
the only independent operations which produce origami numbers. 

\proclaim Theorem. 
$\FF_0 = \FF_{\sqrt{1+x^2}}$. 

\demo Proof. 
Since we already know that $\FF_{\sqrt{1+x^2}} \subset\FF_0$, we only need 
to show that $\FF_0 \subset \FF_{\sqrt{1+x^2}}$. 
That is, we need to show that any origami number may be expressed using 
the usual field operations and the operation $x\mapsto \sqrt{1+x^2}$. 
It is enough to consider the coordinates of origami constructible points, 
because a number is an origami number if and only if it is a coordinate 
of a constructible point. 
There are only four distinct ways of constructing new origami points from 
old ones using the axioms for origami construction.  
These are illustrated in figures~7 and 8. 

\hskip -4bp\psfig{file=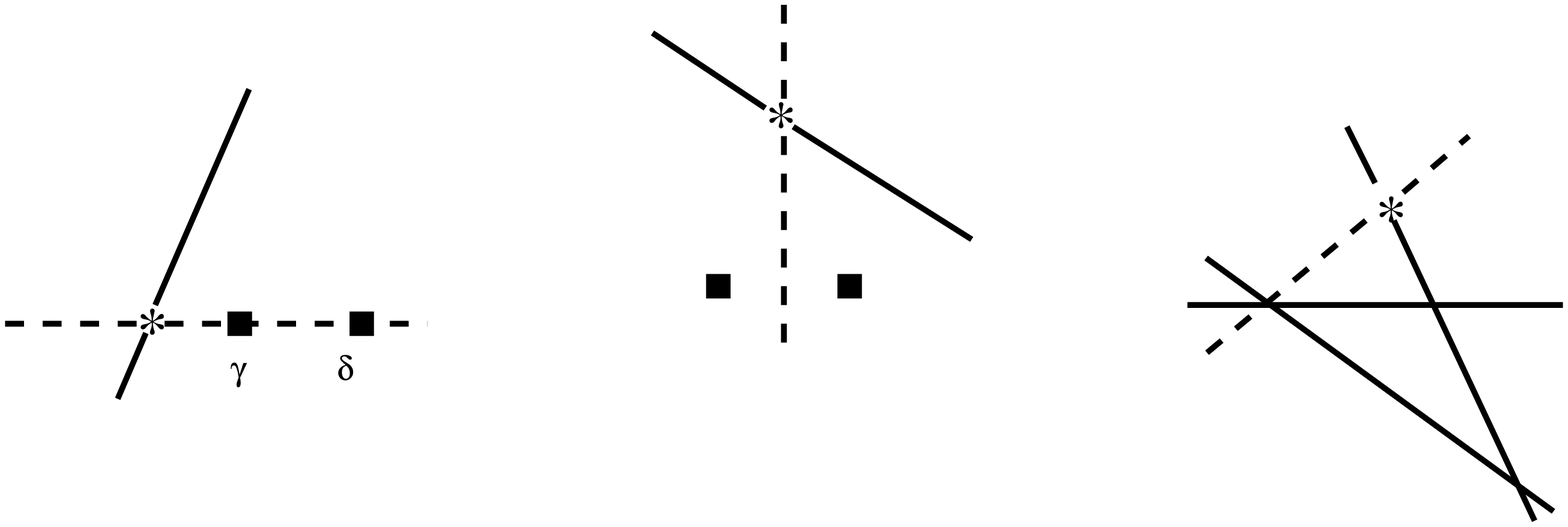,width=5.9truein}

\medskip
\centerline{Figure 7}
\smallskip

\hskip 60bp\psfig{file=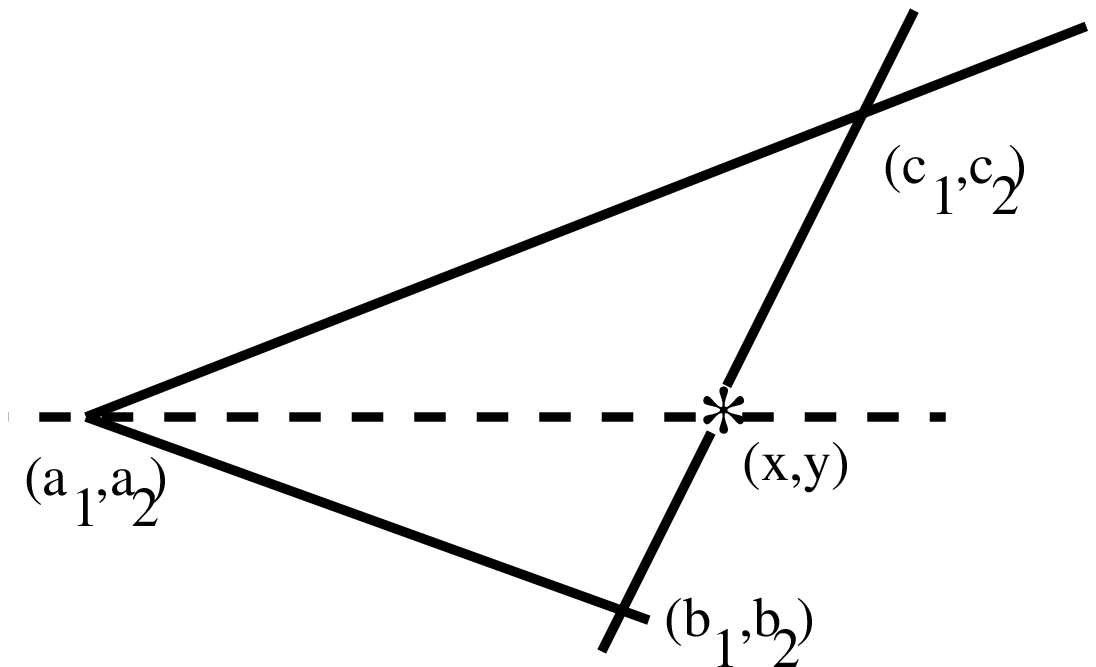,width=2.8truein}

\medskip
\centerline{Figure 8}
\smallskip
\noindent
The only way a new  point will be constructed is by a new crease 
intersecting an old one. 
The four ways of making a crease are:
folding a line between two existing points as in the line $\ov{\gamma\delta}$ 
in figure~7, 
folding the perpendicular bisector to two points as in the second part 
of figure~7, 
reflecting a line as in the third part of figure~7, 
or forming the angle bisector as in figure~8. 
We will explain the case illustrated in figure~8 and leave the remaining 
three cases to the reader. 
When showing that the point $(x,y)$ only depends on the prescribed 
operations, we may assume that $(a_1,a_2) = (0,0)$ by translation, 
because the point $(x,y)$ is found by adding $(a_1,a_2)$ to the 
translated point. 
We may further assume that $(b_1,b_2)$ is on the unit circle, by  
scaling because multiplying by $\sqrt{b_1^2+b_2^2} =|b_1|\sqrt{1+(b_2/b_1)^2}$ 
will reverse the scaling. 
Even further $(b_1,b_2)$ may be assumed to be (1,0) because the 
rotation $(x,y)\mapsto (b_1x-b_2y,b_2x+b_1y)$ sends the point (1,0) 
back to $(b_1,b_2)$. 
Let $\theta = \angle cab$, with the above assumptions, $\cot \theta = c_1/c_2$ 
and $\csc\theta = \sqrt{c_1^2+c_2^2}/c_2 = \sqrt{1+(c_1/c_2)^2}$. 
Now the slope of the new crease is $m=\tan (\theta/2) = \csc\theta 
- \cot \theta$, which only depends on the prescribed operations. 
The new point $(x,y)$ is the intersection of the two lines $y=mx$ and 
$y= [c_2/(c_1-1)](x-1)$, so 
$$x= {c_2\over c_2-m(c_1-1)} \quad \hbox{and}\quad 
y= {mc_2\over c_2-m(c_1-1)}$$ 
which only depends on the prescribed operations as was to be shown.\qed 

The preceding theorem gives an algebraic description of the field 
of origami numbers, and in principle answers which shapes are 
constructible and which are not constructible using origami. 
In practice it is still difficult to decide whether or not a given 
number is an origami number. 
For example, $\sqrt{4+2\sqrt2}$ is an origami number because 
$\sqrt{4+2\sqrt2} = \sqrt{1+(1+\sqrt2\,)^2}$, but what about 
$\sqrt{1+\sqrt2}$? In order to answer this question we need a better 
characterization of origami numbers. 
Before we proceed we will review some elementary facts from abstract algebra 
\ref{AH}, \ref{L}. 

\demo Definition. 
A number, $\alpha$, is an {\it algebraic number\/} if it is a root of a 
polynomial with rational coefficients. 
\vskip1pt
Any algebraic number, $\alpha$, is a root of a unique monic irreducible 
polynomial in $\que [x]$, denoted by $p_\alpha (x)$. 
This polynomial, moreover, divides any polynomial in $\que [x]$ having 
$\alpha$ as a root. 

\demo Definition. 
The {\it conjugates of\/} $\alpha$ are the roots of the polynomial 
$p_\alpha (x)$. An algebraic number is {\it totally real\/} if all of 
its conjugates are real. 
We denote the set of totally real numbers by $\FF_{TR}$. 

Of the numbers which we are using to motivate this section, 
$\sqrt{4+2\sqrt2}$ is totally real, because all of its conjugates 
$(\pm \sqrt{4\pm 2\sqrt2}\,)$ are real, but $\sqrt{1+\sqrt2}$ is not 
totally real because two of its conjugates are imaginary 
$(\pm \sqrt{1-\sqrt2}\,)$. 

The last topic which we review is symmetric polynomials. The symmetric  
group on $n$ letters acts on polynomials in $n$ variables by 
$\sigma f(x_1,x_2,\ldots,x_n) = f(x_{\sigma(1)} ,x_{\sigma(2)},\ldots, 
x_{\sigma(n)})$ where $f\in R[x_1,\ldots,x_n]$ and $R$ is an arbitrary ring. 

\demo Definition. 
The fixed points of the above action are called {\it symmetric polynomials\/}  
over $R$. 

For example, $x_1^2 +x_2^2$ is a symmetric polynomial in two variables 
because it remains unchanged when the variables are interchanged. 
However, $x_1^2 - x_2^2$ is not a symmetric polynomial because it 
becomes $x_2^2 - x_1^2 \ne x_1^2 -x_2^2$ when $x_1$ and $x_2$ are 
interchanged. 
One important class of symmetric polynomials is the class of elementary 
symmetric polynomials. 

\demo Definition. 
If $\prod_{k=1}^n (t+x_k)$ is expanded, we obtain 
$$\prod_{k=1}^n (t+x_k) = \sum_{\ell=0}^n \sigma_\ell (x_1,\ldots,x_n) 
t^{n-\ell}\ .$$ 
The $\sigma_\ell (x_1,\ldots,x_n)$ are the  {\it elementary symmetric 
polynomials\/}. 
\vskip1pt
It is easily verify that 
\vskip1pt
\iitem{} $\sigma_1 = x_1+x_2+\cdots +x_n $
\iitem{} $\sigma_\ell = $ the sum of all products of $\ell$ distinct $x_k$'s 
\iitem{} $\sigma_n = x_1\cdot x_2\cdots x_n$. 

\demo Fact. \ref{L, page 191}. 
The algebra of symmetric polynomials over $R$ is generated by the elementary 
symmetric polynomials. That is, any symmetric polynomial is a linear 
combination of products of the elementary symmetric polynomials. 

We will now begin the final characterization of the origami numbers. 
It happens that all origami numbers are totally real. 
To prove this, it is necessary to show that the sum, difference, product 
and quotient of totally real numbers is totally real, and that 
$\sqrt{1+\alpha^2}$ is totally real whenever $\alpha$ is totally real. 
This is proven by using symmetric polynomials and the following lemma. 

\proclaim Lemma. 
$$\prod_{i=1}^n\ \prod_{j=1}^m (t-x_iy_j) = \det (tI -AB)$$ 
where $A$ and $B$ are matrices with entries expressed in terms of the 
elementary symmetric polynomials of $x_i$ or $y_j$ respectively. 

This lemma is interesting because it is easier to prove a more general 
statement which implies the lemma than it is to verify the lemma. 
We will prove the lemma when the $x_i$ and $y_j$ are independent variables, 
a more general statement than when the $x_i$ and $y_j$ represent numbers, 
but, nevertheless, an easier statement to prove. 

\demo Proof. 
Let 
$$P_A(t) = \prod_{k=1}^n (t-x_k) = \sum_{\ell=0}^n (-1)^\ell \sigma_\ell 
({\bf x}) t^{n-\ell}$$ 
and 
$$P_B(t) = \prod_{j=1}^n (t-y_j) =\sum_{j=0}^m (-1)^j \sigma_j 
({\bf y}) t^{m-j}\ .$$ 
Let 
$$V_{k,\ell} = \left[\matrix{
1\cr 
x_k\cr 
x_k^2\cr 
\vdots\cr 
x_k^{n-1}\cr
y_\ell\cr 
x_k y_\ell\cr  
\vdots\cr 
x_k^{n-1} y_\ell\cr 
\vdots \cr
y_\ell^{m-1}\cr 
x_k y_\ell^{n-1}\cr 
\vdots\cr 
x_k^{n-1} y_\ell^{n-1}\cr}\right]
\quad \raise.75truein\vbox{$
\eqalign{
&\hbox{and let $\ov{A}$ be the $n\times n$ matrix,}\cr
\noalign{\vskip20pt}
&\ov{A} = \left[\matrix{
0&1&0&0&0&\cdots&0\cr 
0&0&1&0&0&\cdots&0\cr 
0&0&0&1&0&\cdots&0\cr 
\vdots&\vdots&\vdots&\vdots&\vdots&&\vdots\cr
(-1)^{n+1}\sigma_n({\bf x})
&(-1)^n \sigma_{n-1}({\bf x})&&&&\cdots&\sigma_1({\bf x})\cr}\right]\cr}
$}
$$
Now let $A$ be the following $nm\times nm$ matrix 
$$A=\left[ \matrix{\ov{A}&&&&\cr &\ov{A}&&&\cr &&\ov{A}&&\cr &&&\ddots&\cr
&&&&\ov{A}\cr}\right]\ .$$ 
By plugging $x_k$ into $P_A(t)$, we find that 
$$x_k^n = \sum_{\ell=1}^n (-1)^{\ell+1} \sigma_\ell ({\bf x}) x_k^{n-1}\ .$$ 
This implies that 
$$AV_{k,\ell} = x_k\cdot V_{k,\ell}$$ 
where $A$ is independent of $k$ and $\ell$. 
In a similar way we can construct a matrix, $B$, with entries given by 
the elementary symmetric functions such that 
$$BV_{k,\ell} = y_\ell V_{k,\ell}\ .$$ 
Now 
$$\eqalign{AB V_{k,\ell} & = Ay_\ell V_{k,\ell}\cr 
&= y_\ell AV_{k,\ell}\cr 
&= x_k y_\ell V_{k,\ell}\ .\cr}$$
Thus $\{ x_k y_\ell\}$ are $nm$ distinct roots of $\det (tI-AB)$ 
which is a monic polynomial of degree $nm$. 
Therefore,  
$$\det (tI-AB) = \prod_{i=1}^n\ \prod_{j=1}^m 
(t-x_i  y_j)\ .$$ 
If the $x_k$'s and $y_\ell$'s were not independent variables, we would not be 
able to conclude that the elements in $\{x_k y_\ell\}$ are distinct.\qed  

With this lemma, we are ready to prove that the set of totally real numbers 
form a field under the operation $x\mapsto \sqrt{1+x^2}$. 

\proclaim Theorem. 
$\FF_{\sqrt{1+x^2}} \subset \FF_{TR}$. 

\demo Proof. 
If $\alpha,\beta\in\FF_{TR}$, we must show that  
$-\alpha$, $\alpha^{-1}$, $\sqrt{1+\alpha^2}$, $\alpha+\beta$, 
$\alpha\cdot \beta\in\FF_{TR}$. 
Let $\{\alpha_i\}_{i=1}^n$ be the conjugates of $\alpha$ and 
$\{\beta_j\}_{j=1}^m$ be the conjugates of $\beta$. 
We will prove the theorem by considering the following five polynomials. 
$$\eqalign{
q_{-\alpha} (t) & = \prod_{i=1}  (t+\alpha_i)\ ,\cr
q_{\alpha^{-1}}(t) &=\biggl(\prod_{i=1}^n (t-\alpha_i^{-1})\biggr) 
\biggl(\prod_{i=1}^n \alpha_i\biggr)\ ,\cr
q_{\sqrt{1+\alpha^2}}(t) & = \prod_{i=1}^n (t^2 - 1-\alpha_i^2)\ ,\cr 
q_{\alpha+\beta} (t) &= \prod_{i=1}^n \prod_{j=1}^m (t-\alpha_i-\beta_j)\ ,\cr 
q_{\alpha\beta}(t) &=\prod_{i=1}^n \prod_{j=1}^m (t-\alpha_i\beta_j)\ .\cr}$$ 
The proofs of the first three cases are similar, and the proofs of the last 
two cases are similar, so we will only prove, in detail, the third case 
and the fifth case. 
If we expand $q_{\sqrt{1+\alpha^2}}(t)$, it is clear that the coefficients 
of $t^k$ will be symmetric polynomials in the $\alpha_i$. 
They may, therefore, be expressed as rational polynomials in the elementary 
symmetric polynomials of the $\alpha_i$. 
Since $(-1)^\ell \sigma_\ell ({\balpha})$ are the coefficients of the 
minimal polynomial for $\alpha$ we may conclude that $q_{\sqrt{1+\alpha^2}}(t) 
\in \que [t]$. It is clear that $\sqrt{1+\alpha^2}$ is a root of 
$q_{\sqrt{1+\alpha^2}}(t)$, thus the minimal polynomial of 
$\sqrt{1+\alpha^2}$, $p_{\sqrt{1+\alpha^2}}(t)$, divides 
$q_{\sqrt{1+\alpha^2}} (t)$. 
The fact that $\alpha$ is totally real implies that all of the conjugates,  
$\alpha_i$, are real. 
Thus, $1+\alpha_i^2$ are all real and positive, so $\pm\sqrt{1+\alpha_i^2}$ 
are all real. 
We now conclude that all of the roots of $q_{\sqrt{1+\alpha^2}}(t)$ are 
real, and therefore $\sqrt{1+\alpha^2}$ is totally real. 
\vskip1pt 
For the fifth case, we use the previous lemma to conclude that 
$q_{\alpha\beta}(t) \in \que[t]$. 
Clearly, $\alpha\beta$ is a root of $q_{\alpha\beta}(t)$ and all of 
the roots of $q_{\alpha\beta}(t)$ are real because $\alpha$ and $\beta$ 
are totally real. 
In the other three cases, it is necessary to show that each of the $q$'s 
are polynomials with rational coefficients and only real roots. 
The first two cases may be tackled with the fact that the elementary 
symmetric polynomials generate the algebra of all symmetric polynomials. 
The fourth case may be verified with a lemma analogous to the 
previous lemma stating that 
$$\prod_{i=1}^n \prod_{j=1}^m (t-x_i-y_j) = \det (tI-A-B)\ .\eqno\eop$$ 

This theorem gives us a practical way to decide that certain shapes 
may not be constructed using origami. 
For example, it is not possible using origami, to construct two cubes 
such that the volume of the second cube is twice that of the first cube. 
If this construction were possible, $\root 3\of 2$ would be an origami 
number and would therefore be totally real. 
One, however, finds that the conjugates of $\root 3\of 2$ are $\root 3\of 2
(-\frac12 \pm\frac{\sqrt3}2 i)$ and $\root 3\of2$, but the first two 
are not real, so $\root 3\of 2$ is not an origami number. 

As we have seen before, $\sqrt2 = \sqrt{1+1^2}$ and $\sqrt{4+2\sqrt2} =
\sqrt{1+(1+\sqrt2\,)^2}$ are origami numbers, so 
$\sqrt{2+\sqrt2} = \sqrt2^{-1} \sqrt{4+2\sqrt2}$ is an origami number. 
>From this we see the  following corollary. 

\proclaim Corollary. 
It is not possible to construct a right triangle with arbitrarily given 
hypotenuse and leg using origami. 

\demo Proof. 
If this were possible, it would be possible to construct a right triangle 
with hypotenuse $\sqrt{2+\sqrt2}$ and leg~1, since these are origami 
numbers. 
Any such triangle would have a leg of length 
$\sqrt{1+\sqrt2} = \sqrt{(\sqrt{2+\sqrt2}\,)^2-1^2}$, but this is 
impossible because $\sqrt{1+\sqrt2}$ is not totally real.\qed  

The following corollary is a consequence of the standard algebraic description 
of compass and straight edge constructions and the two previous theorems 
\ref{AH}. 

\proclaim Corollary. 
Every thing which is constructible with origami is constructible with 
a compass and straight edge, but the converse is not true. 

We want to expand on the relationship between compass and straight edge 
constructions and origami constructions. 
To review, compass and straight edge constructions, let $\FF_{\sqrt x}$ be 
the smallest subfield of $\complex$ closed under the operation 
$x\mapsto \sqrt{x}$, then $\FF_{\sqrt x} \cap \real$ is the collection 
of numbers which are constructible with a compass and straight edge. 
>From our work thus far, it is evident that the origami numbers, $\FF_0$,  
are contained in $\FF_{\sqrt x} \cap \FF_{TR}$. 
It is in fact the case that $\FF_0 = \FF_{\sqrt x} \cap \FF_{TR}$. 
This characterization of the origami numbers is related to David Hilbert's 
$17^{th}$ problem. 
At the International Congress of Mathematics at Paris in 1900, Hilbert 
gave a list of 23 problems \ref{B}. 
His $17^{th}$ problem was to show that any rational function which is 
non-negative when evaluated at any rational number is a sum of squares 
of rational functions. 
In 1926, Artin solved Hilbert's $17^{th}$ problem \ref{Ar}. 
The key idea which Artin used  was the notion of totally positive. 
An element of a field is defined to be totally positive if it is positive 
in every order on the field. 
Artin proved that an element is totally positive if and only if it is a sum 
of squares. 
This is the idea which we use to prove the final characterization of the 
origami numbers. 

\demo Fact. \ref{L, page 457}. 
If $K$ is a finite real algebraic extension of $\que$, then an element 
of $K$ is a sum of squares in $K$ if and only if all of its real conjugates 
are positive. 

\proclaim Theorem. 
$\FF_0 = \FF_{\sqrt{1+x^2}} = \FF_{\sqrt x} \cap \FF_{TR}$. 

\demo Proof. 
We have already shown that $\FF_0 =\FF_{\sqrt{1+x^2}}$ and that 
$\FF_0 \subset \FF_{\sqrt x}\cap \FF_{TR}$, so we need to show that 
$\FF_{\sqrt x} \cap \FF_{TR} \subset \FF_{\sqrt{1+x^2}}$. 
If $\alpha\in \FF_{\sqrt x} \cap \FF_{TR}$, then there exists a 
sequence of totally real numbers, $\{\beta_i\}_{i=1}^n$ and a sequence 
of totally real fields $\{K_j\}_{j=0}^{n-1}$ such that $K_0=\que$, 
$K_i = K_{i-1}(\beta_k)$, $\alpha=\beta_n$, and each $\beta_i$ has 
degree~2 over $K_{i-1}$. 
Since $\beta_i$ has degree~2 over $K_{i-1}$, $\beta_i$ is a root of 
a polynomial of the form 
$$x^2 +c_i x+d_i\ ,$$ 
where $c_i,d_i\in K_{i-1}$. 
Therefore, $(\beta_i +c_i /2)^2 = c_i^2/4- d_i$. 
By the proof of the previous theorem, we know that every conjugate 
of $(\beta_i+c_i/2)^2$ is the square of some conjugate of $\beta_i+ c_i/2$. 
Hence, each of the conjugates of $(\beta_i +c_i/2)^2$ are 
positive and $(\beta_i +c_i/2)^2$ is a sum of squares of elements in $K_{i-1}$. 
Say that 
$$(\beta_i + c_i/2)^2 = r_{i,1}^2 + r_{i,2}^2 +\cdots + r_{i,m}^2\ ,$$ 
then, 
$$\beta_i = r_{i,1} \sqrt{ 1+\left[ {r_{i,2}\over r_{i,1}} 
\sqrt{ 1+\left[ {r_{i,3}\over r_{i,2}} \sqrt{\cdots}\right]^2 } 
\ \right]^2 } - {c_i\over 2}$$
and we are done. 
This shows that any totally real number in $\FF_{\sqrt x}$ is an origami 
number.\qed  

Legend has it that the ancient Athenians were faced with a plague. 
In order to remedy the situation, they sent a delegation to the oracle 
of Apollo at Delos. 
This delegation was told to double the volume of the cubical altar to Apollo. 
However, the Athenians doubled the length of each side of the altar, 
thereby creating an altar with eight times the volume rather than twice 
the volume of the original altar. 
Needless to say, the plague only got worse. 
For years, people have tried to double the size of a cube with compass 
and straight edge, and the gods have not smiled upon them. 
We now can see that the gods will not be satisfied with our elementary 
origami either. 

\Refs 

\item{[AH]} 
G. Alexanderson, A. Hilman, 
``A First Undergraduate Course in Abstract Algebra,'' 3rd  edition,  
Wadsworth Publishing Company, 1983. 

\item{[Ar]} 
E. Artin, 
{\it \"Uber die zerlegung definiter funktionen in quadrate}, 
Abh. Math. Sem. Hausischen Univ. 
{\bf 5} (1927), 100--115. 

\item{[B]} 
F. Browder, 
{\it Mathematical developments arising from Hilbert problems}, 
in ``Proceedings of Symposia in Pure Mathematics,'' Volume 28, 
American Mathematical Society, 1976. 

\item{[K]} F. Klein, 
{\it Vortrage \"uber ausgewahlte Fragen der Elementargeometrie}, 
Teubner, 1895.

\item{[L]} 
S. Lang, 
``Algebra,'' 3rd edition, 
Addison-Wesley, 1993. 

\item{[M]} 
J. Montroll, 
``Origami for the Enthusiast,'' Dover, 1979. 

\item{[R]} T. Sundara Row, 
{\it Geometric Exercises in Paper Folding}, 
Dover, 1966.

\end